\documentclass[10pt, a4paper]{article}

\usepackage{latexsym}
\usepackage{a4wide}
\usepackage{amscd}
\usepackage{graphics}
\usepackage{amsmath}
\usepackage{amssymb}
\usepackage[latin1]{inputenc}
\usepackage{mathrsfs}
\input xy
\xyoption{all}
\usepackage{bm}

\usepackage{geometry}

\newcommand{\N}{\mathbb N}
\newcommand{\R}{\mathbb R}

\newcommand{\Z}{\mathbb Z}
\newcommand{\n}[1]{{\bf #1}}
\newcommand{\Id}{\mathrm{Id}}

\newcommand{\gendeg}{\mathrm{Deg}}
\newcommand{\baricentri}{\Sigma}

\newcommand{\finedim}{\hfill $q.e.d.$\\}            

\newcommand{\proof}{{\sl Proof.}\hspace{5pt}}   

\newtheorem{mainthm}{\sc Theorem}           
\newtheorem{thm}{\sc Theorem}[section]      
\newtheorem{cor}[thm]{\sc Corollary}        
\newtheorem{lem}[thm]{\sc Lemma}            
\newtheorem{prop}[thm]{\sc  Proposition}     
\newtheorem{defn}[thm]{\sc Definition}      
\newtheorem{rem}[thm]{\sc Remark}       
\newtheorem{notation}{\sc Notation}    
\begin{document}

\title{Morse theory for a fourth order elliptic equation with exponential
nonlinearity}
\author{Laura Abatangelo, Alessandro Portaluri\thanks{The  author was partially supported by
PRIN {\em Variational Methods and Nonlinear Differential
Equations.}}}

\date{November 13,  2009}
\maketitle





\begin{abstract}
Given a Hilbert space $(\mathcal H, \langle \cdot, \cdot \rangle)$,
$\Lambda\subset (0, +\infty)$ an interval and $K \in C^2(\mathcal H,
\R)$ whose gradient is a compact mapping, we consider the family of
functionals of the type:
\[
I(\lambda, u) = \dfrac12\langle u, u \rangle -\lambda K(u), \ \
(\lambda, u) \in \Lambda \times \mathcal H.
\]
As already observed by many authors, for the functionals we are
dealing with the (PS) condition may fail under just this
assumptions. Nevertheless, by using a recent deformation Lemma
proven by Lucia in \cite{Luc07}, we prove a Poincar\'e-Hopf type
theorem. Moreover by using this result, together with some
quantitative results about the formal set of barycenters, we are
able to establish a direct and geometrically clear degree counting
formula for a nonlinear scalar field equation on a bounded and
smooth $C^\infty$ region of the four dimensional Euclidean space in
the flavor of \cite{Mal08}. We remark that this formula has been
proven with complete different methods in \cite{LinWeiWan} by using
blow-up type estimates.
\end{abstract}
\smallskip
\centerline{{\em Key Words:\/} Scalar field equations, Geometric
PDE's, Morse Theory, Leray-Schauder degree.}
\bigskip \centerline{\bf AMS
subject classification: 35B33, 53A30, 53C21, 58E05.}
\section*{Introduction}\label{sec:intro}
Let $(\mathcal H, \langle \cdot, \cdot \rangle)$ be a Hilbert space
whose associated norm will be denoted by $\| \cdot\|$. Given an
interval $\Lambda$ of $ (0, \infty)$ and $K$ such that
\begin{equation}\label{eq:condsuKintro}
K \in C^{2}(\mathcal H, \R), \qquad \textrm{with}\quad \nabla K:
\mathcal H \to \mathcal H \ \ \textrm{compact},
\end{equation}
let us consider the functionals which are of the form:
\begin{equation}\label{eq:classefunzionaliintro}
I(\lambda, u) = \dfrac12 \langle u, u \rangle - \lambda K(u), \quad
(\lambda, u)\in \Lambda \times \mathcal H.
\end{equation}
We observe that the conditions
\eqref{eq:condsuKintro}-\eqref{eq:classefunzionaliintro} are not
enough to ensure the (PS)-condition which is known to hold only for
bounded sequences. (See, \cite[Lemma 2.3]{Luc07}). Therefore the
classical flow defined by the vector-field $-\nabla_uI(\lambda, u)$
is not suitable to derive a deformation lemma. However, by using a
recent deformation result proven by \cite[Proposition 1.1]{Luc07},
we prove the following result.
\begin{mainthm}\label{thm:poincarehopfgenintro}
Let $I(\lambda, \cdot)$ be a family of functionals satisfying
\eqref{eq:condsuK}-\eqref{eq:classefunzionali} and 
fix $\bar I(\cdot):= I(\bar \lambda, \cdot)$ for some $\bar\lambda
\in \Lambda$. Given $\varepsilon>0$, let  $\Lambda':=[\bar\lambda
-\varepsilon, \bar\lambda+ \varepsilon]$ be a (compact) subset of
$\Lambda$ and consider $a,b \in \R$ ($a<b$), so that all the
critical points $\bar u$ of $I(\lambda, \cdot)$ for $\lambda \in
\Lambda' $ satisfy $\bar I(\bar u) \in (a,b)$. If
\[
\bar I_a^b:= \left\{u \in \mathcal H:a \leq \bar I(u) \leq b
\right\},
\]
and assuming that $\bar I$ has no critical points at the levels
$a,b$, we have
\begin{equation}\label{eq:formulagradocaratteristicagenintro}
\deg_{LS}(\nabla \bar I, \bar I_a^b, 0)= \chi(\bar I^b, \bar I^a).
\end{equation}
where we denoted by $\deg_{LS}$ the Leray-Schauder degree of the
compact vector field $\nabla \bar I$.
\end{mainthm}
Now, let $\Omega \subset \R^4$ be a bounded connected  and open set
having smooth $C^\infty$ boundary, and let us consider the following
boundary value problem
\begin{equation}\label{eq:probelmaellitticointro}
\left\{\begin{array}{ll} \Delta^2 u = \tau \dfrac{h(x)
e^{u}}{\int_\Omega h(x) e^{u}dx}& \textrm{in} \ \ \Omega\\
u=\Delta u=0\quad & \textrm{on}\  \ \partial \Omega
\end{array}\right.
\end{equation}
where $h$ is a $ C^{2,\alpha}$ positive function for $\alpha \in
(0,1)$, and $\tau \in \R^+$. As already observed by many authors,
the importance of this equation is related to its physical meaning.
In fact, it arises in mathematical physics as a mean field equation
of Euler flows or for the description of self-dual condensates of
some Chern-Simons-Higgs model. (See
\cite{Luc05,Luc06,Luc07,DjaMal,Mal07}, for further details). If
$\mathscr H$ denotes the space of all functions of Sobolev class
$H^2(\Omega) \cap H_0^1(\Omega)$ endowed with the equivalent norm
$\|u\|_\mathscr H:=\|\Delta u\|_2$, than problem
\eqref{eq:probelmaellitticointro} has a variational structure and
for each fixed constant $\tau$, the (weak) solutions can be found as
critical points of the functional
\begin{equation}\label{eq:funzionale}
I_\tau(u):= \dfrac12\|u\|^2_\mathscr H
-\tau\log\left(\dfrac{1}{|\Omega|} \int_\Omega h(x) e^u
dx\right)\qquad \forall \, u \in \mathscr H,
\end{equation}
where we denoted by $|\cdot|$ the Lebesgue measure in $\R^4$. The
key analytic fact which we need in order to classify the critical
points of \eqref{eq:funzionale} is a version for higher order
operators of the Moser-Trudinger inequality. As a direct consequence
of this inequality, it follows that the functional
\eqref{eq:funzionale} is coercive for $\tau < 64\pi^2$ and thus it
is possible to find the solutions of
\eqref{eq:probelmaellitticointro}, by using the direct method of the
calculus of variation. If $ \tau > 64 \pi^2$, the functional
$I_\tau$ is unbounded both from below and from above and hence the
solutions have to be found by other methods, for instance as saddle
points, by using some min-max scheme. A general feature of the
problem is a compactness property if $\tau$ is not integer multiple
of $64\pi^2$ as proven by Lin \& Wei in \cite{LinWei}.\\\noindent If
$\tau < 64\pi^2$ or $\tau \in (64 k \pi^2, 64(k+1)\pi^2)$, $k \in
\N$, by elliptic regularity and by taking into account the
compactness result proven in \cite[Theorem 1.2]{LinWei2}, it is
possible to define the Leray-Schauder degree for the boundary value
problem \eqref{eq:probelmaellitticointro}, fixing a large ball
$\mathscr B_R \subset \mathscr H$ centered at 0 and containing all
the solutions. In fact, let us consider the family of compact
operators $T_\tau: \mathscr B_R \to \mathscr H$, defined by
\[
T_\tau(u):= \tau \,\Delta^{-2} \dfrac{h e^u}{\int_\Omega h e^u};
\]
then the Leray-Schauder degree
\[
d_\tau:= \deg_{LS} (I-T_\tau, \mathscr B_R, 0)
\]
is well-defined for $\tau \neq 64k\pi^2$, $k \in \N$.
\begin{notation}
For any two integers $k_1, k_2$, we use the notation
$\binom{k_1}{k_2}$ to denote
\[
\binom{k_1}{k_2}:=\left\{ \begin{array}{ll}\dfrac{
k_1(k_1-1)\dots(k_1-k_2+1)}{k_2!} & \textrm{if}\ k_1>0\\
1 & \textrm{if}\ k_1=0,
\end{array}\right.
\]
and $\n{k}$ to denote the set $\{1, \dots, k\}$.
\end{notation}
By applying Theorem \ref{thm:poincarehopfgen}, together with a
precise homological properties of the formal set of barycenters
obtained in \cite{DjaMal} we can reprove the following result.
\begin{mainthm}\label{thm:teo2}(\cite{LinWeiWan})
For $\tau \in (64k \pi^2, 64(k+1)\pi^2)$, and $k \in \N$, the
Leray-Schauder degree $d_\tau$ of \eqref{eq:probelmaellitticointro}
is given by
\[
d_\tau = \binom{k-\chi(\Omega)}{k},
\]
where $\chi(\Omega)$ denotes the Euler characteristic of the domain
$\Omega$.
\end{mainthm}
As direct consequence if $\chi(\Omega)\leq 0$ then the problem
\eqref{eq:probelmaellitticointro} possesses a solution provided that
$\tau \not= 64k\pi^2$, $k \in \N$.

In the rest of the section we briefly describe the method and the
main ideas of the proof. As already observed for $\tau > 64 \pi^2$,
the functional $I_\tau$ is unbounded both from above and below due
to the so-called {\em bubbling phenomenum\/} which often occurs in
geometric problems. More precisely, for a given point $x \in \Omega$
and for $\lambda
>0$, we consider the following function
\[
\varphi_{\lambda, x}(y)= \dfrac14\log \left( \dfrac{2\lambda}{1+
\lambda^2 dist(y,x)^2}\right)^4
\]
where $dist(\cdot, \cdot)$ denotes the metric distance on $\Omega$.
For large $\lambda$, one has $e^{\varphi_{\lambda,x}}
\rightharpoonup \delta_x$ (the Dirac mass at $x$) and moreover one
can show that $I(\tau, \varphi_{\lambda,x}) \to - \infty$ as
$\lambda \to + \infty$. Similarly, if $\tau > 64 \pi^2$ it is
possible to construct a function $\varphi$ of the above form (near
at each $x_i$) with $e^{\varphi_{\lambda,\sigma} }\rightharpoonup
\sum_{i=1}^kt_i \delta_{x_i}$ and on which $I_\tau$ still attains
large negative values.  A crucial observation, as proven in
\cite{DjaMal}, is that the constant in Moser-Trudinger inequality
can be divided by the number of regions where $e^u$ is supported.
From this argument we see that one is led naturally to consider the
family  of elements $\sigma:=\sum_{i=1}^k t_i \delta_{x_i}$ with
$(x_i)_i\subset \Omega$ and $\sum_{i=1}^k t_i=1$, known in
literature as the {\em formal set of barycenyters of $\Omega$ of
order $k$\/} and introduced for the first time by Bahri \& Coron in
\cite{BahCor88}. Using the functions $\varphi_{\lambda, x}$, is
indeed possible to map (non-trivially) this set into $\mathscr H$ in
such a way that the functional $I_\tau$  on the image is close to
$-\infty$. On the other hand, it is also possible to do the
opposite, namely to map appropriate sublevels of $I_\tau$ into the
formal set of barycenters. The composition of these two maps turns
out to be homotopic to the identity on the formal set of barycenters
(which is not contractible) and therefore they are both
topologically non-trivial. We remark that our method is along the
same line of a recent result proven by Malchiodi in \cite{Mal08},
for a general Paneitz operator on compact four dimensional
Riemannian manifolds without boundary.

\vskip0.2truecm

\noindent {\bf Acknowledgements.\/}  We would like to thank Prof. A.
Malchiodi for suggesting the problem and for many useful discussions
about this research project.


\section{Preliminaries}\label{sec:nfap}

The aim of this section is to recall some abstract results from
degree theory for $\alpha$-contractions, Sard's lemma for Fredholm
maps and to recall some topological and homological properties of
the so-called {\em formal set of barycenters\/}. Our main references
will be \cite{BahCor88,Dei85,DjaMal,KalKar,Mal08}.

\paragraph{The Sard-Smale theorem and Kuratowski non-compactness
measure.} We start this section with the classical Sard-Smale
theorem stated in a form suitable for our purposes. See
\cite[pag.91]{Dei85}.
\begin{thm}(Sard-Smale)\label{thm:sardsmale}
Let $\Gamma$ be an open subset of a Hilbert space $X$. Suppose that
$\mathscr G \in C^1(\Gamma, X)$ is proper when restricted to any
closed bounded subset of $\Gamma$ and that $\nabla \mathscr
G(x)=\Id-K(x)$ where for every $x \in \Gamma$, $K(x)$ is a compact
operator. Then the set of regular values of $\mathscr G$ is dense in
$X$.
\end{thm}
We will apply this result to $X= \mathscr H$ and $\mathscr G =
\nabla I_\tau$. Since both the map $\mathscr G$ and its  Fr\'echet
derivative are of the form $\Id-K$ where $K$ is a  compact operator,
than the assumptions of theorem \ref{thm:sardsmale} are fulfilled.\\
\noindent  Now let $\Gamma$ be an open subset of $X$ and let
$\mathscr F \colon \Gamma \to X$ be a strict $\alpha$-contraction,
meaning that $\alpha(\mathscr F(B))< k \alpha(B)$ for some fixed $k
\in [0,1)$, where $B\subset \Omega$ is a bounded subset and where
$\alpha$ denotes the {\em Kuratowski measure of non-compactness\/}.
If $y \notin (\Id- \mathscr F)(\partial \Omega)$ and $(\Id- \mathscr
F)^{-1}(\{y\})$ is compact, we can define the {\em generalized
degree\/} $\gendeg$, in such a way that if $\Id- \mathscr F$ is a
compact vector field and $\Gamma$ is a bounded subset it enjoys all
the properties of the Leray-Schauder degree.

\paragraph{Formal set of barycenters.} The aim of this paragraph is to
recall some facts about the formal set of barycenters. \\\noindent
Following \cite{BahCor88}, these spaces are defined by
\begin{equation}\label{eq:definoriginaleformalset}
\baricentri_k:=\left\{\sum_{i=1}^n t_i \delta_{x_i}\  | \ \ (x_1,
\dots, x_n) \in \Omega^n, (t_1, \dots t_n) \in \Delta_{n-1}\right\}
\end{equation}
where $\delta_x$ is the true {\em Dirac mass\/} at the point $x$ and
$\Delta_{n-1}$ is the $(n-1)$-simplex of all tuples $(t_1, \dots,
t_n)$ such that $t_i \geq 0$ and $\sum_{i=1}^n t_ 1=1$. We observe
that the set $\baricentri_k$ is provided by the weak convergence of
measures. In order to give a more topological insight on these
spaces, some definitions are in order.

Let us denote by $J_k$ the $k$-fold join of $\Omega$. We recall that
a point $\underline x \in J_k$ is specified by:
\begin{enumerate}
\item[(i)] $k$ real numbers $t_1, \dots t_k$ satisfying $t_i \geq 0$,
$\sum_{i=1}^k t_i =1$, and
\item[(ii)] a point $x_i \in \Omega$ for each $i\in \n{k}$ such that $t_i \not=0$.
\end{enumerate}
Such a point will be denoted by the symbol $\oplus_{i=1}^k t_ix_i$,
where the elements $x_i$ may be chosen arbitrarily or omitted
whenever the corresponding $t_i$ vanishes. Furthermore we will endow
this space with the strongest topology such that the coordinate
functions are continuous. Now, if $\Sigma^k$ denotes the symmetric
group over $k$ elements, we assume that $\Sigma^k$ acts on $J_k$ by
permuting factors, namely
\[
\forall\  \sigma \in \Sigma^k \colon \ \ \sigma(t_1x_1\oplus
\dots\oplus t_kx_k) := (t_{\sigma(1)} x_{\sigma(1)}\oplus \dots
\oplus t_{\sigma(k)} x_{\sigma(k)}).
\]
Thus, the $k$-th {\em symmetric join\/} of $\Omega$, say $SJ_k$ is
defined as the quotient of $J_k$ with respect to $\Sigma^k$.
\begin{defn}(\cite[Definition 5.1]{KalKar})\label{def:bar}
The $k$-th barycenter space $\baricentri_k$ can be defined as the
quotient of the symmetric join $SJ_k$ under the equivalence relation
$\sim$:
\[
t_1 x_1 \oplus t_2 x_1\oplus \dots \oplus t_kx_k \sim (t_1 +t_2)
x_1\oplus \dots \oplus t_kx_k.
\]
\end{defn}
That is a point in $\baricentri_k$ is a formal abelian sum with the
topology that when $t_i=0$ the entry $0 x_i$ is discarded from the
sum, and when $x_i$ moves in coincidence with $x_j$, one identifies
$t_ix_i +t_jx_i$ with $(t_i+t_j)x_i$. It is possible to show that we
have the embeddings
\[
\Omega \hookrightarrow \baricentri_2 \hookrightarrow \dots
\hookrightarrow \baricentri_{k-1} \hookrightarrow \baricentri_k
\]
and each factor is contractible in the next. Let $P$ be the
projection on $\mathscr H$ (i.e. $P\varphi = \varphi-h$ with
$\Delta^2 h=0$ in $\Omega$ and $h=\varphi$ and $\Delta h =\Delta
\varphi$ on $\partial \Omega$), $\Sigma\subset \mathscr H$ be the
unit sphere and finally let
\[
R: \mathscr H \backslash \{0\} \to \Sigma: u \mapsto
R(u):=u/\|u\|_\mathscr H.
\]
Let $g_k:SJ_k \to \Sigma$ be the map defined as: $ g_k((x_1, \dots,
x_k), (t_1, \dots, t_k)):= R\big(\sum_{i=1}^k t_i P
\varphi_{\lambda, x_i}\big),$ where $\lambda>0$ is fixed and
$\varphi_{\lambda, x_i}$ are given by
\begin{equation}\label{eq:bubblingfunctions}
\varphi_{\lambda, x}(y)= \dfrac14\log \left( \dfrac{2\lambda}{1+
\lambda^2 dist(y,x)^2}\right)^4.
\end{equation}
We observe that since two elements in $SJ_k$ equivalent for the
relation introduced in definition \ref{def:bar} have the same image
through $g_k$, this implies that $g_k$ is well-defined on the
quotient. Denoting by $ \Omega^k$ the $k$-fold product of copies of
$\Omega$ and by $\Delta_k$ the {\em collision set\/} $
\bigcup_{i,j=1}^k \Delta_{i,j}$, where
\[
\Delta_{i,j}:=\{(x_1, \dots, x_k) \in \Omega^k| \ \ x_i=x_j, i
\not=j, \ \ \textrm{for}\ \ i, j \in \n{k}\},
\]
we define the {\em configuration space\/} $ \widehat{\mathfrak
X}_k:= \mathfrak X_k\backslash \Delta_k$. Let us consider the
fibration
\[
\mu \colon \widehat{\mathfrak X}_k \to  \widehat{\mathfrak X}_{k-1},
\quad \textrm{defined by}\ \ \ \mu(x_1, \dots, x_k):=(x_1, \dots,
x_{k-1}),
\]
it is easy to observe that each fiber $\mu^{-1}((x_1, \dots,
x_{k-1}))= \Omega \backslash \{x_1, \dots, x_{k-1}\}$ is
homeomorphic to each other. Thus by using the classical Hopf theorem
for fibrations (see, for instance Spanier \cite{Spa66}, for further
details), the Euler characteristic of $\widehat{\mathfrak X}_k$ can
be computed through the fiber $\Omega \backslash \{x_1, \dots,
x_{k-1}\}$ and $\widehat{\mathfrak X}_{k-1}$. By an easy
calculations it follows that
\begin{equation}\label{eq:careulerspazioconfig}
\chi(\widehat{\mathfrak X}_k)=\chi(\Omega) (\chi(\Omega)-1)\dots
(\chi(\Omega)-k+1).
\end{equation}
\begin{lem}\label{thm:dimensionespaziobaricentri}(well-known)
 The set $\baricentri_k\backslash \baricentri_{k-1}$ is an open smooth manifold of
dimension $5k-1$ for each $k \in \N$.
\end{lem}
\proof The case $k=1$ is trivial, since $\baricentri_1=\Omega$ and
$\Omega$ is a four dimensional manifold being an open subset of
$\R^4$. For $k \geq 2$ the join $J_k$ is a smooth manifold. Since
the action of the symmetric group on $J_k$ is free of fixed points
than the symmetric join is a smooth manifold. Moreover, since
$\baricentri_{k-1}$ is the boundary of $\baricentri_k$, than
$\baricentri_k\backslash \baricentri_{k-1}$ is a smooth open
manifold in which the elements in $\baricentri_k\backslash
\baricentri_{k-1}$ are smoothly parameterized by $4k$ coordinates
locating the points $x_i$ and $k-1$ coordinates identifying the
numbers $t_i$'s. The conclusion immediately follows. \finedim
\begin{lem}(well-known) For any $k \geq 1$, the set
$\baricentri_k$ is a is non contractible stratified set.
\end{lem}
\proof (Sketch). It can be proved by arguing as follows. The case
$k=1$ is trivial. For $k \geq 2$ even if the set $\baricentri_{k-2}$
is not a smooth manifold (actually it is a stratified set) however
it is an ENR which implies that there exists a non trivial (mod 2)
orientation class with respect to its boundary. However by using the
$\mathrm{\check C}$ech-cohomology, and by taking into account that
it is isomorphic to the singular cohomology and over $\Z_2$ to the
singular homology, the thesis follows by using the exactness of the
pair once it is proven that
\[
H_{5k-1}(\baricentri_k;\Z_2) \simeq
H_{5k-1}(\baricentri_k\baricentri_{k-1};\Z_2).
\]
(See, for instance, \cite[Lemma 3.7]{DjaMal}, for further
details).\finedim By using the same arguments as in
\cite[Proposition 5.1]{Mal08}, it can be proven the following
result.
\begin{lem}\label{thm:prelimcarcomp}
Let $\eta>0$ be smaller than the injectivity radius of $\Omega$ and
let $G \colon (0,+\infty) \to (0,+\infty)$ be the non-increasing
function satisfying:
\[
G(t)= \dfrac{1}{t} \ \ \textrm{for} \ \ t \in (0, \eta] \qquad
G(t)=\dfrac{1}{2\eta} \ \ \textrm{for}\ \  t >2\eta.
\]
If $d(x_i, x_j):=dist(x_i, x_j)$ and  $F^*: \baricentri_{k}\setminus
\baricentri_{k-1} \to \R$ as follows
\begin{equation}\label{eq:defdiF}
F^*(\sum_{i=1}^k t_i \delta_{x_i}) := - \sum_{i\not=j} G(d(x_i,
x_j)) - \sum_{i=1}^k \frac{1}{t_i(1-t_i)}.
\end{equation}
Then we have
\begin{equation}\label{eq:formulacaratteristicaEP}
\sum_{i=1}^{5k-1} c_i = (-1)^{k-1}\dfrac{\chi(\widehat{\mathfrak
X}_k)}{k!}
\end{equation}
where $c_i$ denotes the number of critical points of index $i$.
\end{lem}
The following result will be crucial in order to compute the
Leray-Schauder degree of our result.
\begin{prop}\label{thm:calcoloeulerbaricentri}
For any natural number $k$ we have:
\[
\chi(\baricentri_k) = 1- \binom{k-\chi(\Omega)}{k}.
\]
\end{prop}
\proof The proof is given by induction over $k$. The case $k=1$ is
trivial being $\baricentri_1$ homeomorphic to $\Omega$. For $k>1$ we
consider the pair $(\baricentri_k, \baricentri_{k-1})$ and we remark
that the Euler characteristic is additive. Thus
$\chi(\baricentri_k)= \chi(\baricentri_k,\baricentri_{k-1}) + \chi
(\baricentri_{k-1})$. \\
{\em Claim 1. The following formula holds for any natural number $k$
\begin{equation}\label{eq:calcoloeulerocoppia}
\chi(\baricentri_k,\baricentri_{k-1})=
(-1)^{k-1}\binom{\chi(\Omega)-k}{k}.
\end{equation}\/}
Once this is done the proposition easily follows. By Lemma
\ref{thm:dimensionespaziobaricentri} the space
$\baricentri_k\backslash\baricentri_{k-1}$ is an open manifold of
dimension $5k-1$ with boundary $\baricentri_{k-1}$ and by the
definition of $F^*$, Palais-Smale condition holds.\\\noindent
Observe that $\baricentri_{k-1}$ is a deformation retract of the
sublevel $F^*_{-L}:=\{F^* \leq -L\} \cup \baricentri_k$ for $L$
sufficiently large and positive (simply by taking the limit for $L
\to + \infty$). Thus denoting by $\widehat F^* \colon \{F^* \geq
-L\} \to \R$ a non-degenerate function $C^2$-close to the
restriction of $F^*$ to the subset $\{F^* \geq -L\}$, by excision of
the sublevel $F^*_{-L}:=\{F^* < -L\}$ and by the classical
Poincar\'e-Hopf theorem it holds
\[
\chi(\baricentri_k,\baricentri_{k-1})= \sum_{i=1}^{5k-1} (-1)^i c_i.
\]
The thesis follows by formula \eqref{eq:formulacaratteristicaEP} and
\eqref{eq:careulerspazioconfig}.\finedim
\paragraph{Improved Moser-Trudinger inequality.}
Let $C_c^\infty(\Omega)$ be the set of all smooth functions with
compact support in $\Omega$, and let $\mathscr H$ be the completion
with respect to the norm $\|u\|_\mathscr H:=\|\Delta u \|_2$. The
space $\mathscr H$ is a Hilbert space with respect to the scalar
product $ \langle u, v \rangle_\mathscr H := \int_\Omega \Delta u \,
\Delta v dx$ for all $u, v \in \mathscr H,$ and, by the local
regularity theorem and by the Poincar\'e inequality, it follows that
$\mathscr H$ agrees with the space of all functions on $\Omega$ of
Sobolev class $H^2(\Omega) \cap H_0^1(\Omega)$. As immediate
consequence of \cite[Theorem 1.2]{LinWei2} and the Schauder
estimates, the following crucial compactness results hold.
\begin{prop}(\cite[Theorem 1.2]{LinWei2})\label{thm:compattezzaLINWEI}
Let $h \colon \Omega \to \R$ be a positive $C^{2,\alpha}$ function
and $\tau \not= 64k\pi^2$ for $k \in \N$. Then the solutions of
\eqref{eq:probelmaellitticointro} are bounded in
$C^{4,\alpha}(\Omega)$ for any $\alpha \in (0,1)$.
\end{prop}
\begin{prop}(\cite[Lemma 2.1]{LinWei})\label{thm:noboundarybubbles}
Let $u$ be a solution of  \eqref{eq:probelmaellitticointro} with
$\tau \leq c$, for some constant $c$. Then there exists a $\delta>0$
such that
\[
u(x)\leq c,\qquad \forall\, x \,\in \Omega_\delta,
\]
where $\Omega_\delta:=\{x\in \Omega: \  d(x,
\partial \Omega) \leq \delta\}$.
\end{prop}
We remark that proposition \ref{thm:noboundarybubbles} exclude
boundary bubbles.
\begin{lem}\label{thm:lemmaImprovedTrudingerMoser}
There exists a constant $C_\Omega$ depending only on $\Omega$ such
that for all $u \in \mathscr H$ one has:
\begin{equation}\label{eq:mosertrudingerclassica}
\log\left(\dfrac{1}{|\Omega|}\int_\Omega e^{(u- \bar u)}dx
\right)\leq C_\Omega + \dfrac{1}{128 \pi^2}\|u\|^2_\mathscr H
\end{equation}
where $\bar u := \frac{1}{|\Omega|}\int_\Omega u dx$ stands for the
average of $u$ over $\Omega$.
\end{lem}
\proof In fact by \cite[Theorem 1]{Ada88}, there exists $C'_\Omega
>0$ depending only on $\Omega$ such that for all $u \in C_c^2(\Omega)$ it holds
\begin{equation*}\label{eq:giusta}
\dfrac{1}{|\Omega|}\int_\Omega e^{\frac{32\pi^2(u-\bar
u)^2}{\|u\|^2_\mathscr H}}dx \leq C'_\Omega, \qquad \forall \, u \in
\mathscr H.
\end{equation*}
Since for every $a, b \in \R$, we have $(8\pi a -
\frac{1}{8\pi}b)^2\geq 0$ is $ 2ab \leq \frac{1}{64\pi^2}b^2 +
64\pi^2 a^2,$ by setting $a:= u-\bar u$ and $b =\|u\|^2_\mathscr H$,
and exponentiating, we have
\[
\dfrac{1}{|\Omega|}\int_\Omega e^{(u-\bar u)} dx \leq
e^{\frac{1}{128\pi^2}\|u\|^2_\mathscr
H}\dfrac{1}{|\Omega|}\int_\Omega e^{\frac{32\pi^2(u-\bar
u)^2}{\|u\|^2_\mathscr H}} dx \ \leq\
e^{\frac{1}{128\pi^2}\|u\|^2_\mathscr H} C'_\Omega, \qquad \forall
\, u \in \mathscr H.
\]
Taking the logarithm of this last inequality the conclusion follows
by setting $C_\Omega := \log \, C_\Omega $. \finedim In order to
study how the function $e^{u}$ is {\em spread\/} over $\Omega$ we
need some quantitative results. In fact, we will show that if $e^u$
has integral bounded from below on $(l+1)$-regions, the constant
$\frac{1}{128\pi^2}$, can be basically divided by $(l+1)$. The proof
of the proposition \ref{thm:improvedMT},  is up to minor
modifications, an adaptation of the arguments given in \cite[Lemma
2.2]{DjaMal}; we will reproduce it for the sake of completeness.
\begin{prop}\label{thm:improvedMT} For any fixed integer $l$, let
$\Omega_1, \dots, \Omega_{l+1}$ be subsets of $\Omega$ satisfying
$dist(\Omega_i, \Omega_j)\geq \delta_0$, for $i \not=j$, when
$\delta_0$ be positive real number, and let $\gamma_0 \in (0,
\frac{1}{l+1})$. Then for any $\tilde \varepsilon >0$ there exists a
constant $\widetilde C:=\widetilde C(l, \tilde \varepsilon,
\delta_0, \gamma_0)$ such that
\begin{equation*}
\log\left(\dfrac{1}{|\Omega|}\int_\Omega e^{(u-\bar u)}dx\right)\leq
\dfrac{1}{128(l+1)\pi^2-\tilde\varepsilon}\|u\|^2_\mathscr H +
\widetilde C, \qquad
\end{equation*}
for all functions $u \in \mathscr H$ satisfying
\begin{equation}\label{eq:49}
\dfrac{\int_{\Omega_i} e^u \,dx}{\int_\Omega e^u\,dx} \geq \gamma_0
\qquad \forall \, i  \in \n{l+1}.
\end{equation}
\end{prop}
\proof We consider the $(l+1)$ smooth cut-off functions $g_1, \dots,
g_{l+1}$, satisfying the following properties:
\begin{equation}\label{eq:17maldjad}
\left\{\begin{array}{llll} g_i(x) \in [0,1] & \textrm{for every} \ \
x
\in \Omega;\\
g_i(x)=1 & \textrm{for every}\ \  x
\in \Omega_i, i \in \n{l+1};\\
g_i(x)=0 & \textrm{if} \ \ dist(x, \Omega_i) \geq \frac{\delta_0}{4};\\
\|g_i\|_{C^4(\Omega)} \leq C_{\delta_0},
\end{array}\right.
\end{equation}
where $C_{\delta_0}$ depends only on $\delta_0$. By interpolation,
(see, for instance, \cite[Prop. 4.1]{Lio}), for any $\varepsilon
>0$, there exists $C_{\varepsilon, \delta_0}$, such that for any $ v
\in \mathscr H$, and for any $i \in \n{l+1}$ there holds
\begin{equation}\label{eq:18}
\|g_i v\|^2_\mathscr H:=\int_\Omega |\Delta(g_i v)|^2 dx \leq
\int_\Omega g_i^2 |\Delta v|^2 dx + \varepsilon \int_\Omega |\Delta
v|^2 dx + C_{\varepsilon, \delta_0} \int_\Omega v^2 dx.
\end{equation}
Let $u-\bar u = u_1 + u_2$ with $u_1 \in L^{\infty}(\Omega)$, then
from our assumptions we deduce
\[
\int_{\Omega_i} e^{u_2} dx \geq e^{-\|u_1\|_{L^\infty(\Omega)}}
\int_{\Omega_i} e^{(u-\bar u)} dx \geq
e^{-\|u_1\|_{L^\infty(\Omega)}} \gamma_0 \int_\Omega e^{(u-\bar u)}
dx \ \ i \in \n{l+1}.
\]
By invoking inequality \eqref{eq:mosertrudingerclassica} in lemma
\ref{thm:lemmaImprovedTrudingerMoser}, together with the last two
inequalities, it follows that
\begin{eqnarray}
\log\left(\dfrac{1}{|\Omega|}\int_\Omega e^{(u-\bar u)}dx\right)
&\leq& \log\dfrac{1}{\gamma_0} + \|u_1\|_{L^\infty(\Omega)} + \log
\left(\dfrac{1}{|\Omega|}\int_\Omega e^{g_i u_2}dx\right)+
C_\Omega\\\nonumber &\leq& \log\dfrac{1}{\gamma_0} +
\|u_1\|_{L^\infty(\Omega)} +\dfrac{1}{128 \pi^2}\|g_i u\|^2_\mathscr
H + C_\Omega.
\end{eqnarray}
We choose $i$ such that $\int_\Omega |\Delta(g_iu_2)|^2dx \leq
\int_\Omega |\Delta(g_ju_2)|^2dx $, for each $j \in \n{l+1}$. Since
the functions $g_j$ have disjoint supports, the last formula and
\eqref{eq:18}, implies that
\begin{eqnarray*}
\log\left(\dfrac{1}{|\Omega|}\int_\Omega e^{(u-\bar u)}dx\right)
&\leq& \log\dfrac{1}{\gamma_0} + \|u_1\|_{L^\infty(\Omega)} +
C_\Omega+ \left(\dfrac{1}{128(l+1) \pi^2}+\varepsilon\right)\|
u_2\|^2_\mathscr H + \\ &+&C_{\varepsilon, \delta_0} \int_\Omega v^2
dx .
\end{eqnarray*}
Now let $\lambda_{\varepsilon, \delta_0}$ to be an eigenvalue of
$-\Delta^2$ such that $\frac{C_{\varepsilon,
\delta_0}}{\lambda_{\varepsilon, \delta_0}}< \varepsilon$, and we
set
\[
u_1 := P_{V_{\varepsilon, \delta_0}} (u- \bar u); \qquad u_2 :=
P_{V^\perp_{\varepsilon, \delta_0}} (u- \bar u).
\]
Here $V_{\varepsilon, \delta_0}$ is the direct sum of the
eigenspaces of $-\Delta^2$ with Navier boundary conditions and
having eigenvalues less or equal than $\lambda_{\varepsilon,
\delta_0}$ and $P_{V_{\varepsilon, \delta_0}}$,
$P_{V^\perp_{\varepsilon, \delta_0}} $ the orthogonal projections
onto $V_{\varepsilon, \delta_0}$  and $V^\perp_{\varepsilon,
\delta_0}$, respectively. Since $V_{\varepsilon, \delta_0}$ is
finite dimensional, the $L^2$ norm and $L^\infty$ norm of $u - \bar
u$ on $V_{\varepsilon, \delta_0}$ are equivalent; then, by using the
Poincaré-Wirtinger inequality (cfr. \cite[pag. 308]{Brezis}), there
holds:
\[
\|u_1\|^2_{L^\infty(\Omega)} \leq \hat C_{\varepsilon,
\delta_0}\|u_1\|^2_{L^2(\Omega)}\leq \hat C_{\varepsilon,
\delta_0}\|u_1\|^2_{(H^2\cap H^1_0)(\Omega)} \leq \hat
C'_{\varepsilon, \delta_0}\|u_1\|^2_\mathscr H ,
\]
where $C'_{\varepsilon, \delta_0}$ is another constant depending
only on $\varepsilon$ and $\delta_0$. Furthermore
\[
C_{\varepsilon, \delta_0}\int_\Omega u_2^2 dx \leq
\dfrac{C_{\varepsilon, \delta_0}}{\lambda_{\varepsilon,
\delta_0}}\|u_2\|^2_{H^2(\Omega)\cap H^1_0(\Omega)}\leq \varepsilon
\|u_2\|^2_{H^2(\Omega)\cap H^1_0(\Omega)}\leq \varepsilon C'_\Omega
\|u_2\|^2_{\mathscr H},
\]
where $C'_\Omega$ is a constant depending only on $\Omega$. Hence
the last formulas imply
\begin{eqnarray*}
\log\left(\dfrac{1}{|\Omega|}\int_\Omega e^{(u-\bar u)}dx\right)
&\leq& \log\dfrac{1}{\gamma_0} +C'_{\varepsilon,
\delta_0}\|u_1\|_\mathscr H +  C_\Omega+ \left(\dfrac{1}{128(l+1)
\pi^2}+\varepsilon\right)\| u_2\|^2_\mathscr H +\varepsilon
C'_\Omega
\|u_2\|_{\mathscr H}\\
&\leq&\log\dfrac{1}{\gamma_0} + C_\Omega+\left(\dfrac{1}{128(l+1)
\pi^2}+3\varepsilon\right)\| u\|^2_\mathscr H+ \overline
C_{\varepsilon, \delta_0}
\end{eqnarray*}
where $\overline C_{\varepsilon, \delta_0}$ depends only on
$\varepsilon$ and $\delta_0$ (and $l$ which is fixed). This conclude
the proof. \finedim In the next Lemma we show a criterion which
implies the situation described in the first condition in
\eqref{eq:49}.
\begin{lem}(\cite[Lemma 2.3]{DjaMal})\label{thm:covering}
Let $l$ be a given positive integer, and suppose that $\varepsilon$
and $r$ are positive numbers. Suppose that for a non-negative
function $f\in L^1(\Omega)$ with $\|f\|_1=1$ there hold
\[
\int_{\cup_{i=1}^l B_r(p_i)} f dx < 1- \varepsilon, \qquad \forall\
l-\textrm{tuple}\ \  p_1, \dots , p_l \in \Omega.
\]
Then there exists $\bar \varepsilon
>0$ and $\bar r
>0$, depending on $\varepsilon, r, l$ and $\Omega$ (but not on $f$),
and $l+1$ points $\bar p_1, \dots, \bar p_{l+1} \in \Omega$
satisfying
\[
\int_{B_{\bar r}(\bar p_1)}f dx \geq \bar \varepsilon,\dots ,
\int_{B_{\bar r}(\bar p_{l+1})}f dx \geq \bar \varepsilon; \qquad
 B_{2\bar r}(\bar p_i)\cap B_{2\bar r}(\bar
p_j)=\emptyset \ \ for\ i \not=j.
\]
\end{lem}
\begin{lem}\label{thm:Lemma2.4}
If $\tau \in (64k\pi^2, 64(k+1)\pi^2)$ with $k \geq 1$, the
following property holds. For any  $\varepsilon >0$ and any $r>0$
there exists a large positive $L=L(\varepsilon , r)$ such that for
every $u \in \mathscr H$ with $\frac{1}{|\Omega|}\int_\Omega e^u
dx=1$ and $I_\tau( u) \leq -L$ there exist $k$ points $p_{1, u},
\dots, p_{k,u}\in \Omega$ such that
\[
\dfrac{1}{|\Omega|}\int_{\Omega\backslash\cup_{i=1}^k B_r(p_{i, u})}
e^udx <\varepsilon.
\]
\end{lem}
\proof To prove the thesis, we argue by contradiction. Thus, there
exist $\varepsilon, r>0$ and a sequence $(u_n)_n \in \mathscr H$
with $1/|\Omega|\int_\Omega e^{u_n}dx=1$ and $ I_\tau( u_n)\to
-\infty$ such that for every $k$-tuple $p_1, \dots, p_k$ in $\Omega$
we have $1/|\Omega|\int_{\Omega\backslash\cup_{i=1}^k B_r(p_i, u)}
e^udx \geq\varepsilon$. Now applying Lemma \ref{thm:covering} with
$l=k$, $f= e^{u_n}$ and finally with $\delta_0=2\bar r$,
$\Omega_j=B_{\bar r}(\bar p_{j})$ and $\bar \gamma_0= \bar
\varepsilon$ for $j \in \n{k}$ and where the symbols
$\delta_0,\Omega_j,\bar \gamma_0$ were defined in Lemma
\ref{thm:lemmaImprovedTrudingerMoser} and $\bar r,B_{\bar r}(\bar
p_{j}), \bar \varepsilon,  (\bar p_j)_j$ were defined in Lemma
\ref{thm:covering}. By this it follows that, for any given $\tilde
\varepsilon >0$ there exists a constant $C>0$ depending on
$\varepsilon, \tilde \varepsilon$ and on $r$ such that
\begin{equation}\label{eq:bounddasotto}
I_\tau(u_n) \geq \dfrac12\|u_n\|^2_\mathscr H - C \tau
-\dfrac{\tau}{64(k+1)\pi^2-\tilde\varepsilon}\dfrac12\|u\|^2_\mathscr
H,
\end{equation}
where the  constant $ C$ does not depends on $n$. Now since $\tau <
64(k+1)\pi^2$, we can choose $\tilde \varepsilon
>0$ small enough that the number $1-\frac{\tau}{64(k+1)\pi^2-\tilde
\varepsilon }:=\delta'>0$. Therefore the inequality
\eqref{eq:bounddasotto} reduces to
\[
I_\tau( u_n) \geq \dfrac{\delta'}{2}\|u_n\|^2_\mathscr H - C
\tau\geq -K,
\]
where $K$ is a positive constant independent of $n$. This violates
our contradiction assumption, and conclude the proof.\finedim Given
a non-negative $L^1$ function $f$ on $\Omega$, we define the
distance of $f$ from $\baricentri_k$ as
\[
dist(f, \baricentri_k):= \sup\left\{\left|\int_\Omega f\psi dx-
\langle \sigma, \psi\rangle\right|\colon \sigma \in \baricentri_k, \
\ \textrm{and} \ \ \|\psi\|_{C^1(\Omega)}\leq1 \right\},
\]
where we denoted by $\langle \cdot, \cdot \rangle$ the usual duality
product. We also define the set
\[
\mathscr D_{\varepsilon, k}=\{f \in L^1(\Omega) \colon f \geq 0, \
\|f\|_{L^1(\Omega)}=1,\ d(f, \baricentri_k) <\varepsilon\}.
\]
With this notation in mind, by Lemma \ref{thm:Lemma2.4} we deduce
the following.
\begin{lem}\label{thm:4.4}
Suppose $\tau \in (64k\pi^2, 64(k+1)\pi^2)$ with $k \geq 1$. Then
for any $\varepsilon>0$ there exists a large positive
$L=L(\varepsilon)$ such that for all $u \in \mathscr H$ with
$I(\tau,u)\leq -L$ and $1/|\Omega|\int_\Omega e^{u}dx=1$, we have
$dist(e^{u}, \baricentri_k)<\varepsilon$.
\end{lem}
We remark that as a direct consequence of \cite[Theorem
1.2,(ii)]{LinWei2}, the blow-up points $p_{j,u}$ at which the
local-mass is concentrated cannot lie on the boundary of $\Omega$.\\

\section{A topological argument}
The aim of this section is to show that an image of the
$\baricentri_k$ can be mapped into very negative sublevels of the
Euler functional $I_\tau$. Moreover this map is non-trivial in the
sense that it carries some homology. The goal of this section is to
sketch the proof of the following result which is given along the
lines of \cite{Mal07}.
\begin{prop}\label{prop:proposizione2.3}
For any $k \in \N$ and $\tau \in (64k\pi^2, 64(k+1)\pi^2)$, there
exists $L>0$ such that the sublevel $\mathscr H^{-L}$ has the same
homology as $\baricentri_k.$
\end{prop}
The proof of the Proposition \ref{prop:proposizione2.3} will follows
from the homotopy invariance of the homology groups once the
following facts will be established.
\paragraph{\em Mapping $\baricentri_k$ into very low sublevels of
$I_\tau$.\/}  To do so, for $\eta>0$ small enough, consider the
smooth non-decreasing cut-off function $\chi_\eta\colon \R^+ \to \R$
satisfying the following properties:
\begin{equation}
\left\{ \begin{array}{ll} \chi_\eta(t)=t, & \textrm{for}\ \ t \in
[0, \eta];\\
\chi_\eta(t)=2\eta, & \textrm{for}\ \ t \geq 2\eta;\\
\chi_\eta(t) \in [\eta, 2 \eta],&  \textrm{for}\ \ t \in [\eta,
2\eta].
\end{array}\right.
\end{equation}
Then given $\sigma \in \baricentri_k$, and $k>0$, we define the
family  of maps $\phi_\lambda \colon \baricentri_k \to H^2(\Omega)$
as $\phi_\lambda(\sigma):=\varphi_{\lambda, \sigma}(\cdot)$ where
the function $\varphi_{\lambda, \sigma} \colon \Omega' \to \R$ is
defined by
\begin{equation}\label{eq:testfunctions}
\varphi_{\lambda, \sigma}(y) :=\dfrac14\log \sum_{i=1}^k t_i
\left(\dfrac{2\lambda}{1+\lambda^2\chi^2_\eta(d_i(y))}\right)^4,
\qquad y \in \Omega
\end{equation}
where we set $ d_i(y)=d(y, x_i)$, for $x_i \in \Omega$. a
Let
\begin{equation}\label{eq:testgen}
\overline \varphi_{\lambda, \sigma}(y):=P(\varphi_{\lambda,
\sigma})(y), \qquad \forall \ y \in \Omega.
\end{equation}
where $P$ is the projection defined above.
We observe that, since the distance function is a 1-Lipschitz
function than $\varphi_{\lambda, \sigma}$ is also a Lipschitz
function in $y$ and hence it belongs to $H^2(\Omega)$.
\begin{prop}\label{thm:mappareibaricentrineisottolivelli}
Let $\overline \varphi_{\lambda, \sigma}$ be defined as in
\eqref{eq:testgen}. Then as $\lambda \to + \infty$ the following
properties hold
\begin{enumerate}
\item[(i)] $e^{\overline \varphi_{\lambda, \sigma}}\rightharpoonup \sigma$ weakly in
the sense of distributions;
\item[(ii)] $I_\tau(\overline \varphi_{\lambda, \sigma})\to -\infty$ uniformly
with respect to  $\sigma \in \baricentri_k$.
\end{enumerate}
\end{prop}
\proof To prove $(i)$ we first consider the function
\[
 \varphi_{\lambda, \sigma}(y):=
\left(\dfrac{2\lambda}{1+\lambda^2\chi^2_\eta(d_i(y))}\right)^4,
\qquad \forall \ y \in \Omega,
\]
where $x$ is a fixed point in $\Omega$. It is easy to verify that
$\varphi_{\lambda, \sigma}(y)\to \delta_{x_i}$ for $\lambda \to +
\infty$. Then $(i)$ follows from the explicit expression of
$\varphi_{\lambda, \sigma}$.\\ \noindent In order to prove (ii), we
evaluate separately each term of $I_\tau$, and claim that the
following estimates hold
\begin{equation}\label{eq:41mal08}
\log\left(\dfrac{1}{|\Omega|}\int_\Omega e^{\overline
\varphi_{\lambda, \sigma}}dx\right)=O(1) \qquad \textrm{as} \quad
\lambda \to + \infty.
\end{equation}
\begin{equation}\label{eq:lemma4.2mal08}
\dfrac12\|\overline \varphi_{\lambda, \sigma}\|_{\mathscr H}^2 \leq
(128k\pi^2 +o_\epsilon(1))\log \lambda + C_{\epsilon}) \qquad
(\textrm{uniformly in $\sigma \in \Sigma_k$}),
\end{equation}
where $o_\epsilon(1)\to 0$ as $\epsilon \to 0$ and where
$C_{\epsilon}$ is a constant independent $(x_i)_i$.\\
\noindent The proof of \eqref{eq:41mal08} it is easy and it follows
by integrating over $\Omega$. The proof of \eqref{eq:lemma4.2mal08}
is much more involved and it follows by Lemma 4.2 in \cite{DjaMal}.
\finedim
\paragraph{\em Mapping very low sublevels of $I_\tau$ into
$\baricentri_k$ and an homotopy inverse.\/}  The main idea is to
construct a non-trivial continuous map $\psi\colon \mathscr H \to
\baricentri_k$ from the sublevels of the Euler functional into
$\baricentri_k$ such that the composition $ \psi \circ \phi_\lambda$
is homotopic to identity on ${\baricentri_k}$.
\begin{prop}\label{thm:mapparesottolivellineibaricentri}
Suppose that $\tau \in (64k\pi^2, 64(k+1)\pi^2)$ with $k \geq 1$.
Then there exists $L>0$ and a continuous projection $\psi\colon
\mathscr H^{-L}\to \baricentri_k$ with the following properties.
\begin{enumerate}
\item[(i)] If $(u_n)_n\subset\mathscr H^{-L}$ is such that
$e^{u_n}\rightharpoonup\sigma$, for some $\sigma \in \baricentri_k$,
then $\psi(u_n) \to \sigma$;
\item[(ii)] if $\overline \varphi_{\lambda, \sigma}$ is as in
\eqref{eq:testgen}, then for any $\lambda$ sufficiently large the
map $\sigma \mapsto \psi(\overline\varphi_{\lambda,\sigma})$ is
homotopic to the identity on $\baricentri_k$.
\end{enumerate}
\end{prop}
\proof First of all we observe that item $(i)$ follows directly from
item $(ii)$ and Proposition
\ref{thm:mappareibaricentrineisottolivelli}. The non-trivial part is
the construction of the global continuous projection map $\psi$
which is a left homotopy inverse has proven in \cite[Section
3]{DjaMal}. \finedim We close this section by observing that, up to
minor modifications, the above defined map $\psi$ is also a right
inverse homotopy as proven in \cite[Appendix]{Mal08}. Thus summing
up we conclude that
\begin{cor}
Given $L$ sufficiently large the topological spaces $\mathscr
H^{-L}$ and $\baricentri_k$ are equivalent, up to homotopy.
\end{cor}

\section{A Poincar\'e-Hopf Theorem without (PS)}

The aim of this section is to prove an analogous of the
Poincar\'e-Hopf theorem for a special class of functionals. To do
so, let $(\mathcal H, \langle \cdot, \cdot \rangle)$ be a Hilbert
space whose associated norm will be denoted by $\| \cdot\|$. Given
an interval $\Lambda$ of $ (0, \infty)$ and a map $K$ such that
\begin{equation}\label{eq:condsuK}
K \in C^{2}(\mathcal H, \R), \qquad \textrm{with}\quad \nabla K:
\mathcal H \to \mathcal H \ \ \textrm{compact},
\end{equation}
let us consider the functionals which are of the form:
\begin{equation}\label{eq:classefunzionali}
I(\lambda, u) = \dfrac12 \langle u, u \rangle - \lambda K(u), \quad
(\lambda, u)\in \Lambda \times \mathcal H.
\end{equation}
It is well-known (see, for instance, \cite[Lemma 2.3]{Luc07}) that
the conditions \eqref{eq:condsuK}-\eqref{eq:classefunzionali} could
de not enough to ensure the (PS)-condition which is known to hold
only for {\em bounded sequences\/}. Now by using the deformation
Lemma proven in \cite[Proposition 1.1]{Luc07}, we are in position to
derive the following result.
\begin{thm}\label{thm:poincarehopfgen}{\em (A Poincar\'e-Hopf theorem)\/}
Let $I(\lambda, \cdot)$ be a family of functionals satisfying
\eqref{eq:condsuK}-\eqref{eq:classefunzionali} and 
fix $\bar I(\cdot):= I(\bar \lambda, \cdot)$ for some $\bar\lambda
\in \Lambda$. Given $\varepsilon>0$, let  $\Lambda':=[\bar\lambda
-\varepsilon, \bar\lambda+ \varepsilon]$ be a (compact) subset of
$\Lambda$ and consider $a,b \in \R$ ($a<b$), so that all the
critical points $\bar u$ of $I(\lambda, \cdot)$ for
$\lambda \in \Lambda' $ satisfy $\bar I(\bar u) \in (a,b)$. 
Assuming that $\bar I$ has no critical points at the levels $a,b$,
we have
\begin{equation}\label{eq:formulagradocaratteristicagen}
\deg_{LS}(\nabla \bar I, \bar I_a^b, 0)= \chi(\bar I^b, \bar I^a).
\end{equation}
\end{thm}
The proof of this result will be given into two main steps. In the
first step we will assume that all the critical points are
non-degenerate; in the second step we will remove this assumption.\\

\noindent \proof {\em First step: non-degenerate case.\/} We let
$\mathscr K$ denote the set of critical points of $\bar I$ which is
compact by hypothesis. By compactness and non-degeneracy
assumptions, $\bar I$ has only finitely-many critical levels each of
whose consists only of finitely-many critical points. Let $R$ be so
large that all the critical points of $I_\lambda$ for $\lambda \in
\Lambda'$ are in $\mathscr B_{\frac{R}{2}}(0)$.  Then we can define
the cut-off function $\theta : \mathcal H \to [0,1]$ satisfying
\[
\theta(u)=1 \ \textrm{for}\ u \in \mathcal B_R(0); \qquad
\theta(u)=0 \ \textrm{for}\ u \in \mathcal H \backslash \mathcal
B_{2R}(0).
\]
Following Lucia in \cite{Luc07}, let $Z \in C^1(\mathcal H, \mathcal
H)$ be defined by:
\[
Z(u):= -[|\nabla K(u)|\nabla \bar I(u) + |\nabla \bar I(u)|\nabla
K(u)],
\]
and choose $\omega_\varepsilon \in C^\infty(\R)$ such that
\[
0 \leq \omega_\varepsilon \leq 1, \qquad \omega_\varepsilon(\zeta)=0
\ \textrm{for all}\ \ \zeta \leq \varepsilon, \qquad
\omega_\varepsilon(\zeta)=1 \ \textrm{for all} \ \ \zeta \geq
2\varepsilon.
\]
Finally we can define
\[
W(u):=-\omega_\varepsilon\left(\dfrac{|\nabla \bar I(u)|}{|\nabla
K(u)|}\right)\nabla \bar I(u) + Z(u),
\]
where $\omega_\varepsilon\left(|\nabla \bar I(u)|/|\nabla
K(u)|\right)$ is understood to be equal $1$ when $\nabla K(u)=0$.
Given the vector field:
\[
\widetilde W(u):= -\theta(u) \nabla \bar I(u) + (1- \theta(u)) W(u),
\]
we observe that it decreases $\bar I$ in the complement of $\mathscr
K$. We consider the local flow $\eta=\eta(t, u_0)$ defined by the
Cauchy problem:
\[
\dfrac{du}{dt} = \widetilde W(u), \qquad  u(0)=u_0.
\]
{\em Claim 1.If $\bar I$ has no critical levels inside some interval
$[\tilde a, \tilde b]$, then the sub-level $\bar I^{\tilde a}$  is a
deformation retract of $\bar I^{\tilde b}$.} \\ \noindent To prove
this, we arguing as follows. Given $u_0 \in \bar I^{\tilde b}$, we
can prove that
\begin{equation}\label{eq:complicata} \bar I(\eta(t, u_0)) \leq -c^2t + \bar I
(u_0).\footnote{%
The proof of this inequality is the most involved part of this claim
and it can be proven up to minor modifications reapeating word by
word the arguments given in \cite[pagg. 121-122]{Luc07}.\/}
\end{equation} Thus there exists a $t$ such that $\bar
I(\eta(t, u_0) \leq \tilde a$. Then we define:
\begin{equation*}
t_a(u_0):=\left\{\begin{array}{ll}\inf\{t \geq0: \bar I(\eta(t,
u_0)) \in \bar I^{\tilde a}\}& \textrm{if}\quad  \bar I(u_0) > \tilde a\\
0 & \textrm{if}\quad  \bar I(u_0) \leq \tilde a.
\end{array}\right.
\end{equation*}
Now the map
\[
\widetilde \eta :[0,1] \times \mathcal H \to \mathcal H, \qquad (s,
u_0)\mapsto \eta(st_{\tilde a}(u_0), u_0),
\]
is a deformation retraction, as required.\\ \noindent
Now let $\bar c_i$ be the number of critical points of $\bar I$ of
index $i$. By classical Morse-theoretical arguments as in
\cite[Theorem 3.2, 3.3, pagg. 100-103]{Cha93}, by excising $\{\bar I
< \tilde a\}$, it follows that
\[
\deg_{LS}(\nabla \bar I, \bar I_a^b, 0)= \sum_i (-1)^i \bar c_i =
\chi(\bar I^b, \bar I^a).
\]
This conclude the proof in the non degenerate case.\\
\noindent{\em Second step: degenerate case.\/} We reduce ourselves
to the non-degenerate case. To do so, fix a small $\delta >0$ so
that $dist(\mathscr K, \bar I_a^b)> 4\delta$, and define the set
$\mathscr K_\delta = \{u \in \mathcal H: dist(u, \mathscr
K)<\delta\}$. We next choose a smooth cut-off function $p$ such that
\[
p(u)= 1 \ \ \textrm{for every}\  u \in \mathscr K_\delta; \qquad
p(u)=0 \ \ \textrm{for every}\  u \in \mathcal H \backslash \mathscr
K_{2\delta}.
\]
We can also choose $p$ such that $0 \leq p(u) \leq 1$ for all $u \in
\mathcal H$ and having uniformly bounded derivative in $\mathscr
K_{2 \delta}$. Now let $\mathscr G:= \nabla \bar I\vert_{\mathscr
K_\delta}: \mathscr K_\delta \to \mathcal H$. Since the map
$\mathscr G$ is a compact perturbation of the identity, by applying
the Sard-Smale theorem (see theorem \ref{thm:sardsmale}), the set of
regular values of $\mathscr G$ is dense in $\mathcal H$. This
implies that we can find an arbitrarily small $u_0$ such that
$\nabla \mathscr G(p)$ is non-degenerate for each $p \in \mathscr
G^{-1}(u_0)$ which is equivalent to say that $\nabla^2 \bar I$ is
non-degenerate on the set
\[
\Gamma(u_0):=\{u \in \mathcal H: \nabla I(u)=u_0\}\cap \mathscr
K_\delta.
\]
Moreover we observe that $\|\nabla \bar I\|\geq \gamma_\delta>0$ on
$\mathscr K_{2\delta}\backslash \mathscr K_\delta$ for some constant
$\gamma_\delta$. Now let us consider the function
\[
\widetilde I(u):= \bar I(u) + p(u)\langle u_0, u\rangle.
\]
It can be shown that the following facts hold:
\begin{enumerate}
\item[(i)] $\widetilde I$ coincides with $\bar I$ in
$\mathcal H\backslash \mathscr K_{2\delta}$;
\item[(ii)] $\widetilde I$ has the same critical points as $I(\tau, \cdot)$ in
$ \mathscr H\backslash \mathscr K_\delta$;
\item[(iii)] $\widetilde I$ is non-degenerate in $\bar I_a^b$.
\end{enumerate}
Item (i) is trivial. To prove (ii) we observe that since $\widetilde
I$ and $\bar I$ coincides out of $\mathscr K_{2\delta}$, it is
enough to prove the claim for $u \in \overline{\mathscr
K_{2\delta}}\backslash \mathscr K_\delta$. By differentiating, we
have
\[
\langle \nabla \widetilde I(u), v\rangle = \langle \nabla \bar I(u)+
\nabla p(u)\langle u, u_0\rangle +p(u) u_0, v \rangle, \qquad
\forall v \in \mathcal H.
\]
Thus, by recalling that $u \in \overline{\mathscr
K_{2\delta}}\backslash \mathscr K_\delta$, it follows that
\[
\|\nabla \widetilde I(u)\|\geq \|\nabla \bar I(u)\|- |\langle u,
u_0|\|\nabla p(u)\|-p(u)\|u_0\|\geq \gamma_\delta-\|u_0\|(\|\nabla
p(u)\|\|u\|+1)>0,
\]
where the last inequality follows since $p$ has uniformly bounded
derivatives and $u_0$ can be chosen arbitrarily small. To prove
(iii) we argue as follows. Since all the critical points of
$\widetilde I$ are in $\mathscr K_\delta$, let us assume by
contradiction that $\widetilde I$ is degenerate at some critical
point $\bar u$. Now since $\bar u \notin \mathscr K$, it follows
that $\bar u \in \mathscr K_\delta\backslash \mathscr K$. Moreover
$\nabla \widetilde I(\bar u)=0$ is equivalent to say that $\nabla
\bar I(\bar u)=u_0$ and therefore $\bar u \in \Gamma(u_0)$. But this
is contradict the fact that $\nabla^2 \widetilde I(p)$ is
non-degenerate on $p \in \Gamma(u_0)$.

Now, for $\|u_0\|$ sufficiently small the map $\nabla I - Id$ is a
strict $\alpha$-contraction. (See Section \ref{sec:nfap}) and since
$(\nabla \widetilde I)^{-1}(\{u_0\})=\mathscr K$, the generalized
degree $\gendeg(\nabla \widetilde I, \bar I_a^b, u_0)$ is
well-defined; moreover it coincides with $\gendeg(\nabla \widetilde
I, \bar I_a^b, 0)$ since it is locally constant. With the above
choice for $R$ and by using the excision property and the homotopy
invariance of the generalized degree, (see, for instance,
\cite{Dei85} for further details), we have
\[
\deg_{LS}(\nabla \widetilde I, \mathscr B_R, 0)= \gendeg(\nabla
\widetilde I, \mathscr B_R, 0).
\]
Now choosing a possibly larger $R$ in such a way $\mathscr K_{2
\delta} \subset \mathscr B_{R/2}$, the conclusion readily follows by
the first step, simply by replacing $\bar I$ with $\widetilde I$.
\finedim
\begin{cor}\label{thm:corollario2.8}
If $\tau \in (64k\pi^2, 64(k+1)\pi^2)$ for some $k \in \N$ and if
$b$ is sufficiently large positive, the sub-level $\mathscr H^b$ is
a deformation retract of $\mathscr H$ and hence it has the homology
of a point.
\end{cor}
\proof This result follows, by using the deformation constructed in
the proof of the Poincar\'e-Hopf theorem. See, for instance
\cite[Corollary 2.8]{Mal08}.\finedim Setting
\begin{equation}\label{eq:J}
J(u) := \log \left(\dfrac{1}{|\Omega|}\int_\Omega h(x) e^udx\right)
\end{equation}
the functional \eqref{eq:funzionale} can be put in the following
form: $I_\tau(u)=\dfrac12\|u\|^2_\mathscr H- \tau J(u)$.

\paragraph{Proof of Theorem \ref{thm:teo2}.}
\proof In order to prove \ref{thm:teo2}, it is enough to apply
theorem \ref{thm:poincarehopfgen} to the functional
\eqref{eq:funzionale} for $\lambda = \tau$, $\Lambda = (64k \pi^2,
64(k+1)\pi^2)$ for $k \geq1$, $\mathcal H= \mathscr H$ and finally
$K(u)= J(u)$ where $J$ was given in \eqref{eq:J}. The only thing it
should be noted, is that all the critical points $\bar u$ of
$I_\tau$ for $\tau \in [\bar \tau-\varepsilon, \bar
\tau+\varepsilon] \subset (64k \pi^2, 64(k+1)\pi^2) $ satisfy $\bar
I(\bar u) \in (a,b)$. This is a consequence of proposition
\ref{thm:compattezzaLINWEI} and of the boundedness of $J$ which is
consequence of the Moser-Trudinger inequality. Now the conclusion
follows choosing $a=-L$ as in proposition \ref{prop:proposizione2.3}
and $b$ as in corollary \ref{thm:corollario2.8}. In fact by using
theorem \ref{thm:poincarehopfgen}, we have that
\[
d_\tau=\chi(\bar I^b, \bar I^a)=\chi(\bar I^b)- \chi(\bar
I^a)=\chi(\mathscr H) -\chi(\baricentri_k)=1- \chi(\baricentri_k).
\]
The conclusion follows by invoking proposition
\ref{thm:calcoloeulerbaricentri}. \finedim
\begin{rem}
We observe that the Leray-Schauder degree in the Sobolev space
$\mathscr H$ coincides with the degree in every H\"older space
$C^{2,\alpha}(\Omega)$, $\alpha \in (0,1)$. See for instance
\cite[Part I, Appendix B, Theorem B.1-B.2]{Li}.
\end{rem}


\begin{thebibliography}{99}
\bibitem{Ada88} David R. Adams {\em A sharp inequality of J. Moser
for Higher order derivatives\/}. Annals of Math. {\bf 128}, No.2
(1988) pp 385--398.

\bibitem{Bah84} A. Bahri, {\em Un problème variationnel sans
compacité dans la géometrie de contact\/} C.R. Acad. Sci. Paris.
Sér. I Math {\bf 299}, (1984), No. 15, pp. 757--760.


\bibitem{BahCor88} A. Bahri, J. M. Coron {\em On a nonlinear
elliptic equation involving the critical Sobolev exponent: the
effect of the topology of the domain.\/} Comm. Pure Appl. Math.,
{\bf 41} (1988), pp. 253--294.

\bibitem{Brezis} H. Brezis, {\em Analisi Funzionale. Teoria e
Applicazioni.\/} Serie di Matematica e Fisica. Liguori Editore.

\bibitem{Cha93} K.C.Chang {\em Infinite-dimensional Morse theory and
multiple solution problems\/} PNLDE 6. Birkha\"user Boston, Inc.,
Boston, MA, 1993.

\bibitem{Cla44} C.E. Clark, {\em The symmetric join of a complex.\/}
Bull. Amer. Math. Soc. 50, (1944), pp. 81--88.

\bibitem{Dei85} K. Deimling, {\em Nonlinear function analysis\/}.
Springer-Verlag. Berlin (1985).


\bibitem{DjaMal} Z. Djadli, A. Malchiodi {\em  Existence of conformal metrics with constant
$Q$-curvature\/}. Ann. of Math. (2) 168 (2008), no. 3, 813--858.

\bibitem{KalKar} S. Kallel, R. Karoui, {\em Symmetric joins and
weighted barycenters\/}. ArXiv:math/0602283v2

\bibitem{LinWei} C.Lin, J. Wei, {\em Sharp estimates for bubbling solutions of a
fourth order mean field equation.\/} Ann. Sc. Norm. Super. Pisa Cl.
Sci. (5) 6 (2007), no. 4, 599--630.

\bibitem{Li} Y.Y. Li, {\em Prescribed scalar curvature on $S^n$ and related problems.
Part I.} J. Differential Equations 120 (1995), no. 2, 319--410. Part
II. Existence and compactness. Comm. Pure Appl. Math. 49 (1996), no.
6, 541--597.

\bibitem{LinWei2} C.Lin, J. Wei, {\em Locating the peaks of solutions via the maximum principle. II.
A local version of the method of moving planes.\/} Comm. Pure Appl.
Math. 56 (2003), no. 6, 784--809.

\bibitem{LinWeiWan} C.Lin, J. Wei, L. Wang {\em Topological Degree
of fourth order mean field equations}. Preprint 2007.

\bibitem{Lio} Lions, J.-L., {\em Équations différentielles opérationnelles et problèmes aux
limites.\/} Die Grundlehren der mathematischen Wissenschaften, Bd.
111 Springer-Verlag, Berlin-Göttingen-Heidelberg 1961, 292 pp.


\bibitem{Luc05} M. Lucia, {\em A mountain pass theorem without
Palais-Smale condition\/}, C.R. Acad. Sci. Paris, Ser. I {\bf 341}
(2005), pp. 287--291.

\bibitem{Luc06} M. Lucia, {\em A blowing-up branch of solutions for a mean field equation,\/}
Calc. Var {\bf 26}, (2006),No.3 pp. 313--330.

\bibitem{Luc07} M. Lucia, {\em A deformation Lemma with an application to a mean field equation,\/}
Topol. Methods in Nonlinear Analysis {\bf 30} (2007),N0. 1,  pp.
113--138.

\bibitem{Mal07} A. Malchiodi {\em Topological methods for an elliptic equation with
exponential nonlinearities\/}. Discrete Contin. Dyn. Syst. {\bf
21}(2008), No.1, pp.277--294.

\bibitem{Mal08} A. Malchiodi {\em Morse theory and a scalar field
equation on compact surfaces\/}. Adv. Diff. Eq., 13 (2008),
1109-1129.

\bibitem{Mil56} J. Milnor, {\em Construction of Universal bundle
II.\/} Annals of Math. {\bf 63}, No.3, (1956), pp. 430--436.

\bibitem{Spa66} E. H. Spanier {\em Algebraic topology},
Mc-Graw-Hill, New York-Toronto-London, 1966.

\end{thebibliography}
\end{document}